\documentclass[a4paper, 12pt]{article}
\usepackage[english]{babel}
\usepackage{amsfonts}
\usepackage{amssymb}
\usepackage[leqno]{amsmath}
\usepackage{amsthm}
\usepackage{amscd}
\usepackage{xypic}
\usepackage[curve]{xy}
\usepackage{epsfig}
\usepackage{enumerate}

\input xy
\xyoption{all}
\newtheorem{theorem}{Theorem}[section]
\newtheorem{proposition}[theorem]{Proposition}

\newtheorem{corollary}[theorem]{Corollary}

\newtheorem{example}[theorem]{Example}
\newtheorem{lemma}[theorem]{Lemma}
\newtheorem{remark}[theorem]{Remark}
\newtheorem{convention}[theorem]{Convention}

\def\blfootnote{\xdef\thefnmark{}\footnotetext}

\def\Rad{\mathrm{Rad}}

 \DeclareMathOperator{\Hom}{Hom}
\def\Coker{\mathrm{Coker\hspace{0.1cm}}}
\def\Cohom{\mathrm{Cohom}}

\def\Ext{\mathrm{Ext}}
\def\Hom{\mathrm{Hom}}
\def\Rad{\mathrm{Rad}}

\def\dim{\mathrm{dim}}

\def\limdir{\varinjlim}

\def\End{\mathrm{End}}
\def\Soc{\mathrm{Soc\hspace{0.1cm}}}
\def\Socle{\mathrm{Soc}}

\def\T{\mathcal{T}}

\def\M{\mathcal{M}}

\title{Some remarks on localization in coalgebras\blfootnote{E-mail: gnavarro@ugr.es}}
\author{Gabriel
Navarro\footnote{Research supported by Spanish MEC grant
BES-2002-2403, DGES BMF2001-2823 and FQM-266
 (Junta de Andaluc{\'\i}a Research Group).}
 }

\date{\vspace{-30pt}}

\begin{document}

\maketitle

\begin{center} Department of Algebra,  Faculty of Science, University of
Granada \\ c. Fuentenueva s/n,   18071 Granada \\ Spain \end{center}

\abstract{We analyze the geometry of the Ext-quiver of a coalgebra
$C$ in order to study the behavior of simple and injective
$C$-comodules under the action of the functors associated to a
localizing subcategory of the category of $C$-comodules. }

%\abstract{Given an arbitrary coalgebra $C$, we study the behavior of
%simple and injective right $C$-comodules under the action of the
%functors associated to a localizing subcategory of the category of
%right $C$-comodules.}

\section*{Introduction}

In recent years, the difficulty of studying coalgebras in a general
framework has caused the appearance of several kinds of coalgebras
which are investigated separately. Many of them are defined through
certain properties involving their category of comodules or, merely,
their simple or injective comodules. For instance, this is the case
of cosemisimple,
 pure semisimple, semiperfect, quasi-co-Frobenius, hereditary or serial coalgebras, see
\cite{cuadra-serial}, \cite{dascabook}, \cite{jlms}, \cite{lin} or
\cite{simsonnowak}. Therefore it is natural to ask about the
behavior of such comodules in different situations. In this paper we
shall focus on this question in the context of localization as
defined by Gabriel in \cite{gabriel}. In particular, we will study
the behavior of simple and (indecomposable) injective comodules
under the action of the localizing functors associated to a
localizing subcategory. In many cases this analysis is related to
the geometry of the Ext-quiver associated to the coalgebra. For that
reason, in Section \ref{ext-quiver} we take into consideration some
geometric properties of it. We devote Section \ref{section},
\ref{quotient} and \ref{colocalizing} to study how the quotient,
section and colocalizing functors transform simple and injective
objects, respectively. Lastly, in Section \ref{secsemicentral}, we
relate those properties to the idempotent associated to the
localizing subcategory, extending the results of \cite{jmnr}.

\section{Preliminaries and notation}

Throughout we fix a ground field $K$ and we assume that all vector
spaces are over $K$ and every map is a $K$-linear map. In
particular, $C$ is a $K$-coalgebra which, without loss of
generality, we shall consider basic. The reader is expected to be
familiar with basic coalgebra theory, see for instance \cite{Abe},
\cite{montgomery} or \cite{sweedler}. Unless otherwise stated, all
$C$-comodules will be right $C$-comodules. We denote by $\M^C$ the
category of right $C$-comodules and by $\M^C_{qf}$ and $\M^C_{f}$
the full subcategories of $\M^C$ formed by quasi-finite and finite
dimensional comodules, respectively. We also denote by
$\{S_x\}_{x\in I_C}$ and by $\{E_x\}_{x\in I_C}$ a complete set of
pairwise non-isomorphic simple and indecomposable injective right
$C$-comodules, respectively. For each right $C$-comodule $M$, we may
calculate its socle, $\Soc M$, and its injective envelope, $E(M)$.
Then, for each $x\in I_C$, we assume that $\Soc E_x=S_x$, and
consequently, $E(S_x)=E_x$.

We will follow the quiver notation of \cite{jmn} and \cite{jmnr}, in
particular, by a quiver we mean an oriented graph $Q=(Q_0,Q_1)$,
where $Q_0$ is set of vertices and $Q_1$ is the set of arrows. The
path coalgebra $KQ$ is defined in the standard way. The set of
simple right $KQ$-comodules is $\{Kx\}_{x\in Q_0}$, and, for each
$x\in Q_0$, $E_x=E(Kx)$ is generated by the set of paths ending at
$x$.

Following \cite{gabriel}, a dense subcategory $\T$ of $\M^C$ is said
to be localizing if the \emph{quotient functor} $T:\M^C\rightarrow
\M^C/\T$ has a right adjoint functor, $S$, called the \emph{section
functor}. $T$ is exact, $S$ is fully faithful and left exact, and
$TS=1_{\M^C/\T}$. Dually, see \cite{blas2}, $\T$ is said to be
colocalizing if $T$ has a left adjoint functor, $H$, called the
\emph{colocalizing functor}. $H$ is a fully faithful and right exact
functor such that $TH=1_{\M^C/\T}$. It is well-known that there
exist one-to-one correspondences between localizing subcategories,
indecomposable injective comodules and simple comodules. In
\cite{cuadra}, \cite{jmnr} and \cite{woodcock} localizing
subcategories are described by means of idempotent elements of the
dual algebra $C^*$. In particular, it is proved that the quotient
category is the category of right comodules over the coalgebra
$eCe$, where $e$ is the idempotent element associated to the
localizing subcategory. See the above references for the coalgebra
structure of $eCe$ and a description of the localization functors by
means of $e$. The localization in categories of comodules over path
coalgebras is described in detail in \cite{jmnr}.

\section{The geometry of the Ext-quiver}\label{ext-quiver}

To any coalgebra $C$, we may associate a quiver $\Gamma_C$ known as
the right Ext-quiver of $C$, see \cite{mont}. We recall that the set
of vertices of $\Gamma_C$ is the set of pairwise non-isomorphic
simple right $C$-comodules $\{S_x\}_{x\in I_C}$ and, for two
vertices $S_x$ and $S_y$, there exists a unique arrow
$S_y\rightarrow S_x$ in $\Gamma_C$ if and only if
$\Ext^1_C(S_y,S_x)\neq 0$. We may proceed analogously with left
$C$-comodules and obtain the left Ext-quiver of $C$, $\Delta_C$.
Throughout we assume that $\Gamma_C$ is connected, i.e., $C$ is
indecomposable as coalgebra.

\begin{remark} The results obtained in this paper are also valid for the valued Gabriel
quiver of $C$, $(Q^C,d^C)$, i. e., following \cite{sim-kosak}, the
valued quiver
 whose set of vertices is $\{S_x\}_{x\in I_C}$ and such
that there exists a unique valued arrow $\xymatrix{ S_y
\ar[rr]^-{(d^1_{yx},d^2_{yx})} & & S_x}$ if and only if
$\Ext_C^1(S_y,S_x)\neq 0$ and $d^1_{yx}=\dim_{\End_C(S_y)}
\Ext_C^1(S_y,S_x)$ as right $\End_C(S_y)$-module and
$d^2_{yx}=\dim_{\End_C(S_x)} \Ext_C^1(S_y,S_x)$ as left
$\End_C(S_x)$-module.

This is also true if we consider the (non-valued) Gabriel quiver.
That quiver is obtained taking the same set of vertices and the
number of arrows from a vertex $S_x$ to a vertex $S_y$ is
$\dim_{\End_C(S_x)} \Ext^1_C(S_x,S_y)$. We recall that if $C$ is
pointed (or $K$ is algebraically closed) then it is isomorphic to a
subcoalgebra of the path coalgebra of its Gabriel quiver.
\end{remark}

Let us take into consideration some geometric properties of
$\Gamma_C$. Given a vertex $S_x$, we say that the vertex $S_y$ is an
\emph{immediate predecessor} (or a 1-predecessor) of $S_x$ if there
exists an arrow $S_y\rightarrow S_x$ in $\Gamma_C$.

\begin{lemma}\label{predecessors} $S_y$ is an immediate predecessor of $S_x$ if and only
if $S_y\subseteq \Socle (E_x/S_x)$.
\end{lemma}
\begin{proof}
Let us consider the short exact sequence $S_x\hookrightarrow
E_x\rightarrow E_x/S_x$. We apply to it the functor $\Hom_C(S_y,-)$
and then we obtain the long exact sequence
$$\xymatrix@R=6pt{
0 \ar[r]            & \Hom_C(S_y, S_x) \ar[r]^-f        & \Hom_C(S_y,E_x) \ar[r] & \\
 \ar[r] & \Hom_C(S_y, E_x/S_x) \ar[r] & \Ext_C^1(S_y,S_x) \ar[r]  & 0}$$
Since $f$ is a $K$-linear isomorphism, $\Hom_C(S_y, E_x/S_x)\cong
\Ext_C^1(S_y,S_x)$ and the result follows.
\end{proof}

We may generalize the former definition by means of the socle
filtration. Following \cite{green}, any right $C$-comodule $M$ has a
filtration
$$0\subset \Soc M\subset \Socle^2 M\subset \cdots \subset M$$
called the \emph{Loewy series}, where, for $n>1$, $\Socle^nM$ is the
unique subcomodule of $M$ satisfying that $\Socle^{n-1}M\subset
\Socle^nM$ and
$$\frac{\Socle^{n}M}{\Socle^{n-1}M}=\Socle\left (\frac{M}{\Socle^{n-1}M}\right ).$$

\begin{lemma}\cite{green}\label{lemmasocle} Let $M$ and $N$ be right $C$-comodules then:
\begin{enumerate}[$(a)$]
\item Let $n$ be a positive integer. If $0\subseteq R_1\subseteq R_2 \subseteq \cdots \subseteq
R_n$ is a chain of subcomodules of $M$ such that $R_i/R_{i-1}$ is
semisimple for all $i=1,\ldots ,n$ then $R_n\subseteq \Socle^n M$.
\item If $N$ is a subcomodule of $M$ then $\Socle^n N=\Socle^n M
\cap N$ for all $n\geq 1$.
\item If $f:N\rightarrow M$ is a morphism of right $C$-comodules then
$f(\Socle^n N)\subseteq \Socle^n M$ for all $n\geq 1$.
\item If $M=\bigoplus_{\lambda} M_\lambda$ then $\Socle^n
M=\bigoplus_\lambda \Socle^n M_\lambda$ for all $n\geq 1$.
\end{enumerate}
\end{lemma}

We shall need the following result:

\begin{lemma}\label{quotientsocle} Let $M$ be a right $C$-comodule. Then, for each positive integer $n$,
 the chain
$$0\subset \frac{\Socle^{n+1}M}{\Socle^nM}\subset
\frac{\Socle^{n+2}M}{\Socle^nM}\subset \cdots\subset
\frac{\Socle^{n+t}M}{\Socle^nM}\subset \cdots $$ is the Loewy series
of the right $C$-comodule $M/\Socle^n M$, that is,
$$\Socle^t\left (\frac{M}{\Socle^n M} \right
)=\frac{\Socle^{n+t} M}{\Socle^n M} \hspace{1cm} \text{for each
$t>0$}$$
\end{lemma}
\begin{proof}
For each $t>0$, denote by $N_t$ the $C$-comodule
$\Socle^{n+t}M/\Socle^n M$ and by $N$ the $C$-comodule $M/\Socle^n
M$. The case $t=1$ follows from the definition of the Loewy series.
Assume now that the statement holds for $t-1$. Then $N_t$ is a
subcomodule of $N$ such that $\Socle^{t-1} N\subset N_t$ and
$$\frac{N_t}{\Socle^{t-1} N}\cong \frac{\Socle^{n+t}
M}{\Socle^{n+t-1}M}\cong \Socle \left (\frac{M}{\Socle^{n+t-1} M}
\right )\cong \Socle \left (\frac{N}{\Socle^{t-1} N} \right ).$$
Thus $\Socle^t N=N_t$.
\end{proof}

Given a vertex $S_x$, we say that the vertex $S_y$ is an
$n$-\emph{predecessor} of $S_x$ if $\Ext^1_C(S_y, \Socle^nE_x)\neq
0$, or equivalently, proceeding as in Lemma \ref{predecessors}, if
$S_y\subseteq \Soc(E_x/\Socle^{n}E_x)$.

\begin{convention} Throughout, given a simple comodule $S_x$, for each $n\geq 1$, we shall denote by
 $\{S_i\}_{i\in I_{n}}$ the set of simple $C$-comodules such that
 $$\Socle \left (\frac{E_x}{\Socle^{n}
E_x} \right ) =\bigoplus_{i\in I_n} S_i,$$ and we will refer to it
as the set of all $n$-predecessors of $S_x$. %although there were some
%$i,j\in I_n$ such that $S_i=S_j$.
 Obviously, each
$n$-predecessor $S_i$ is repeated $r_i$ times, where
$$r_i=\frac{\dim_K
\Hom_C(S_i,E_x/\Socle^n E_x)}{\dim_K\End_C(S_i)}=\frac{\dim_K
\Ext^1_C(S_i, \Socle^n E_x)}{\dim_K\End_C(S_i)}.$$
\end{convention}

\begin{lemma}\label{equivpredec} Let $S_x$ and $S_y$ be two simple
$C$-comodules. The following assertions are equivalent:
\begin{enumerate}[$(a)$]
\item $S_y$ is a $n$-predecessor of $S_x$.
\item There exists a non-zero morphism $f:\Socle^{n+1} E_x
\rightarrow E_y$ such that $f(\Socle^n E_x)=0$.
\item There exists a morphism $g:E_x\rightarrow E_y$ such that
$g(\Socle^i E_x)=0$ for all $i=1,\ldots ,n$ and $g(\Socle^{n+1}
E_x)\neq 0$ \end{enumerate}
\end{lemma}
\begin{proof} $(a)\Leftrightarrow (b)$.
Assume that $S_y$ is an $n$-predecessor of $S_x$. Then there exists
a non-zero map $h:\Socle^{n+1}E_x/\Socle^{n}E_x\rightarrow E_y$
making commutative the following diagram
$$\xymatrix{ S_y \ar@<-0.3ex>@{^{(}->}[d]_i \ar@<-0.3ex>@{^{(}->}[r]^-i &
\frac{\Socle^{n+1}E_x}{\Socle^{n}E_x} \ar@{-->}[ld]^h
\\ E_y}$$
Then, the composition $\xymatrix{\Socle^{n+1} E_x\ar[r]^-p &
\Socle^{n+1}E_x/\Socle^{n}E_x \ar[r]^-h & E_y}$ is a nonzero
morphism which vanishes in $\Socle^n E_x$.

Conversely, given such an $f$, it decomposes through a non-zero
morphism $g:\bigoplus_{i\in I_n} S_i \longrightarrow E_y$. Therefore
there is an $i\in I_n$ such that $gu_i:S_i\rightarrow E_y$ is
non-zero, where $u_i$ is the standard inclusion. That is, $S_i=S_y$
is an $n$-predecessor of $S_x$.

$(b)\Leftrightarrow (c)$. Assume that $f$ is such a morphism. Since
$E_y$ is an injective $C$-comodule, there exists a morphism
$g:E_x\rightarrow E_y$ such that
$$\xymatrix{ \Socle^{n+1} E_x \ar@<-0.3ex>[d]_f \ar@<-0.3ex>@{^{(}->}[r]^-i &
E_x \ar@{-->}[ld]^g
\\ E_y}$$
commutes. Obviously, $g(\Socle^n E_x)=0$ and $g(\Socle^{n+1}
E_x)\neq 0$.

For the converse, it is enough to consider the restriction of $g$ to
the subcomodule $\Socle^{n+1}E_x$.
\end{proof}

\begin{remark} Observe that the proof is also valid for an arbitrary
$C$-comodule $M$, that is, the following assertions are equivalent:
\begin{enumerate}[$(a)$]
\item $S_y\subseteq M/\Socle^{n}M$.
\item There exists a non-zero morphism $f:\Socle^{n+1} M
\rightarrow E_y$ such that $f(\Socle^n M)=0$.
\item There exists a morphism $g:M\rightarrow E_y$ such that
$g(\Socle^i M)=0$ for all $i=1,\ldots ,n$ and $g(\Socle^{n+1}M)\neq
0$ \end{enumerate}
\end{remark}

We simply
 say that $S_y$ is a predecessor of $S_x$ if there exists an
integer $n\geq 1$ such that $S_y$ is an $n$-predecessor of $S_x$.
The following result gives a necessary and sufficient condition for
the vertex $S_y$ to be a predecessor of a vertex $S_x$. Given two
indecomposable injective right $C$-comodules $E_x$ and $E_y$, we
denote by $\Rad_C(E_x,E_y)$ the set of all morphisms in
$\Hom_C(E_x,E_y)$ which are not bijective.

\begin{corollary}\label{predechom} Let $S_x$ and $S_y$ be two simple
$C$-comodules. Then, $S_y$ is a predecessor of $S_x$ if and only
$\Rad_C(E_x,E_y)\neq 0$.
\end{corollary}
\begin{proof} The sufficiency is proved by the former lemma.

Conversely, for each $n\geq 1$, we have the short exact sequence
$$\xymatrix{\Socle^{n}E_x \ar[r]& \Socle^{n+1}E_x \ar[r] & \bigoplus_{i\in
I_n} S_i}.$$ If $S_y$ is not a predecessor of $S_x$ then $S_y\ncong
S_i$ for all $i\in I_n$. Therefore, for any $n\geq 1$,
$\Hom_C(\Socle^{n+1}E_x,E_y)\cong \Hom_C(\Socle^{n}E_x,E_y)$. Now,
$$\begin{array}{rl} \Hom_C(E_x,E_y)  \cong&\Hom_C(\limdir
\Socle^nE_x,E_y)\\  \cong&\limdir\Hom_C(\Socle^nE_x,E_y)\\
\cong& \Hom_C(S_x,E_y)\\
\cong& \left\{
      \begin{array}{ll}
        0, & \text{if $S_x\ncong S_y$,}\\
        \End_C(S_x), & \text{if $S_x\cong S_y$.}
      \end{array}
    \right.
\end{array}$$
Thus $\Rad_C(E_x,E_y)=0$.
\end{proof}

\begin{theorem}\label{predecpath} Let $S_x$ and $S_y$ be two simple $C$-comodules and $n$ be a positive integer.
If $S_y$ is an $n$-predecessor of $S_x$ then there exists a path in
$\Gamma_C$ of length $n$ from $S_y$ to $S_x$.
\end{theorem}
\begin{proof}
We proceed by induction on the integer $n$. For $n=1$ is just Lemma
\ref{predecessors}. Let us assume that the assertion holds for $n-1$
and let $S_y$ be a $n$-predecessor of $S_x$. By Lemma
\ref{equivpredec}, there exists a non-zero map $f:\Socle^{n+1} E_x
\rightarrow E_y$ such that $f(\Socle^{n}E_x)=0$. In particular we
may decompose $f$ as follows $$\xymatrix{\Socle^{n+1} E_x\ar[d]^-p
\ar[r]^-{f} & E_y \\ M=\frac{\Socle^{n+1}E_x}{\Socle^{n-1}E_x}
\ar[ru]_-g & }$$ where $g$ is a non-zero map such that
$g(\Socle^{n}E_x /\Socle^{n-1} E_x)=0$.

By Lemma \ref{quotientsocle}, $M=\Socle^2(E_x/\Socle^{n-1}E_x)$ and
then it is contained in $\Socle^2(\bigoplus_{i\in I_{n-1}}
E_i)=\bigoplus_{i\in I_{n-1}} \Socle^2 E_i$.  Therefore, since $E_y$
is injective, there exists a non-zero map $h$ making commutative the
following diagram
$$\xymatrix{E_y &  \\ M \ar[u]_{g} \ar[r]^-i & \bigoplus_{i\in I_{n-1}} \Socle^2 E_i
\ar[ul]_-{h} }$$ Hence there is an index $j\in I_{n-1}$ such that
the composition $hu_j:\Socle^2 E_j\rightarrow E_y$ is non-zero,
where $u_j$ is the usual inclusion. Lastly, since $g(\Soc M)=0$,
$h(\bigoplus_{i\in I_{n-1}} S_i)=0$ and then $hu_j(S_j)=0$. Thus
there is an arrow $S_y\rightarrow S_j$ in $\Gamma_C$.
\end{proof}

\begin{remark} The reader should observe that if there is a path in $\Gamma_C$
from $S_y$ to $S_x$, then $S_y$ does not have to be a predecessor of
$S_x$. For example, consider the quiver $Q$
$$\xy\xymatrix@C=30pt{\circ \ar[r]^-{\alpha} & \circ
\ar[r]^-{\beta} & \circ }
\POS (0,-2)*+{\txt{\scriptsize $y$}}
\POS (14.5,-2)*+{\txt{\scriptsize $z$}}
\POS (29,-2)*+{\txt{\scriptsize $x$}}
\POS (32,0)*+{\txt{\scriptsize ,}}
\endxy$$ and the subcoalgebra $C$
of $KQ$ generated by $\{x, y, z, \alpha , \beta\}$. Then the quiver
$\Gamma_C$ is
$$S_y\longrightarrow S_z\longrightarrow S_x.$$
Obviously, there is a path from $S_y$ to $S_x$, but there is no
non-zero morphisms
$$f:E_x=<x,\beta
>\longrightarrow E_y=<y>.$$

On the other hand, if $C$ is the coalgebra $KQ$, the Ext-quiver of
$KQ$ is also the previous quiver but, in this case, we may obtain a
map
$$f:E_x=<x,\beta, \beta\alpha>\longrightarrow E_y=<y>$$ defined by
$f(\beta\alpha)=y$ and zero otherwise. \end{remark}

\begin{lemma} Let $E_x$ and $E_y$ be two indecomposable injective
right $C$-comodules and $f:E_x\rightarrow E_y$ be a morphism of
$C$-comodules such that $f(\Socle^nE_x)=0$ and
$f(\Socle^{n+1}E_x)\neq 0$. Then $f(\Socle^{n+1}E_x)=S_y$ and
$f(\Socle^{n+t}E_x)\subseteq \Socle^t E_y$ for any $t> 1$. Moreover,
if $C$ is hereditary, then $f(\Socle^{n+t}E_x)=\Socle^t E_y$ for any
$t\geq 1$ .
\end{lemma}
\begin{proof}
We may factorize $f$ as the composition:
$$\displaystyle \xymatrix{ E_x\ar[rr]^-p
\ar@/_15pt/[rrrr]_-{f} & & \displaystyle{ \frac{E_x}{\Socle^{n}E_x}}
\ar[rr]^-g & & E_y}$$ Then $f(\Socle^{n+1} E_x)=gp(\Socle^{n+1}
E_x)=g(\frac{\Socle^{n+1}E_x}{\Socle^nE_x})$ is a non-zero
semisimple subcomodule of $E_y$, i.e., it is $S_y$. Let us consider
the chain
$$0\subseteq f(\Socle^nE_x)\subseteq f(\Socle^{n+1}E_x)\subseteq\cdots
\subseteq f(\Socle^{n+t}E_x) \subseteq\cdots$$ Since each quotient
$\frac{f(\Socle^{n+i}E_x)}{f(\Socle^{n+i-1}E_x)}$ is semisimple for
any $i\geq 0$, by Lemma \ref{lemmasocle},
$f(\Socle^{n+t}E_x)\subseteq \Socle^t E_y$ for any $t>1$.

Suppose now that $C$ is hereditary. Then, for each $i\in I_n$, the
following diagram commutes $$\xymatrix{ E_x \ar[rd]_-p \ar[r]^-f & E_y & \\
&\frac{E_x}{\Socle^{n}E_x}\cong \bigoplus_{i\in I_n} E_i \ar[u]^-g&
E_i \ar[l]^-{u_i} \ar[lu]_-{gu_i}}$$ Since $g(\bigoplus_{i\in
I_n}S_i)=S_y$, there is an index $j\in I_n$ such that
$g_{|S_j}:S_j\rightarrow S_y$ is bijective. Thus $gu_j$ is an
isomorphism and
$$\Socle^t E_y=gu_j(\Socle^t E_j)\subseteq g(\bigoplus_{i\in I_n} \Socle^t E_i)
=f(\Socle^{n+t} E_x)$$ for any $t>0$.
\end{proof}

\begin{corollary}Let $C$ be a hereditary coalgebra and $n$ be a
positive integer. The following conditions are equivalent:
\begin{enumerate}[$(a)$]
\item There is a path in $\Gamma_C$ of length $n$ from a vertex $S_y$
to a vertex $S_x$.
\item  $S_y$ is an $n$-predecessor of $S_x$.
\end{enumerate}
\end{corollary}
\begin{proof} It is enough to prove $(a)\Rightarrow (b)$.
If $$\xymatrix{S_y \ar[r] & S_1 \ar[r] & \cdots \ar[r]& S_{n-1}
\ar[r] & S_x}$$ is a path in $\Gamma_C$, there exists a sequence of
(surjective) morphisms
$$\xymatrix@R=30pt{E_x \ar[r]^-{f_{n}} & E_{n-1} \ar[r]^-{f_{n-1}} & \cdots
\ar[r]^-{f_2} & E_1 \ar[r]^-{f_1} & E_y}$$ such that $f_i(S_i)=0$
and $f_i(\Socle^2 E_i)\neq 0$ for all $i=1 ,\ldots, n$, where
$S_n=S_x$. Then, applying repeatedly the previous lemma, we obtain
that  $(f_1f_2\cdots f_n)(\Socle^n E_x)=0$ and $(f_1f_2\cdots
f_n)(\Socle^{n+1}E_x)=S_y\neq 0$. By Lemma \ref{equivpredec}, $S_y$
is an $n$-predecessor of $S_x$.
\end{proof}

\section{The section functor}\label{section}

From now on we fix an idempotent element $e\in C^*$. We will denote
by $\T_e$ the localizing subcategory associated to $e$ and by
$\{S_x\}_{x\in I_e\subset I_C}$ the subset of simple comodules of
the quotient category. Let us consider the quotient and the section
functor associated to $\T_e$:
\[ \xymatrix@C=50pt{ \mathcal{M}^{C}
\ar@<0.75ex>[rr]^-{T=e(-)=-\square_{C}eC} & &
\ar@<0.75ex>[ll]^-{S=-\square_{eCe}Ce} \mathcal{M}^{eCe}}.\] We
recall that there exists a torsion theory on $\mathcal{M}^{C}$
associated to the functor $T$, where a right $C$-comodule $M$ is a
torsion comodule if $T(M)=0$. If $M$ is not torsion, we denote by
$t(M)$ the torsion subcomodule of $M$.

We know that, by \cite{jmnr}, for a simple right $C$-comodule $S_x$,
$T(S_x)=S_x$ if $x\in I_e$ and zero otherwise. From this fact we
obtain the following result:

\begin{lemma}\label{soclequotient}
Let $M$ be a right $C$-comodule then $T(\Soc M)\subseteq \Soc T(M)$.
\end{lemma}
\begin{proof}
Let us suppose that $\Soc M=(\bigoplus_{i\in I} S_i )\bigoplus
(\bigoplus_{j\in J} T_j )$, where $S_i$ and $T_j$ are simple right
$C$-comodules such that $T(S_i)=S_i$ and $T(T_j)=0$ for all $i\in I$
and $j\in J$. Since $\Soc M\subseteq M$ then we have that $T(\Soc
M)=\bigoplus_{i\in I} S_i \subseteq T(M)$.
\end{proof}

Let us study the behavior of the injective comodules under the
action of the section functor. Indeed, we shall prove that $S$
preserves indecomposable injective comodules and, consequently,
injective envelopes. In what follows we will denote by
$\{\overline{E}_x\}_{x\in I_e}$ a complete set of pairwise
non-isomorphic indecomposable injective right $eCe$-comodules, and
assume that $\overline{E}_x$ is the injective envelope of the simple
right $eCe$-comodule $S_x$ for each $x\in I_e$.

\begin{proposition}\label{injectivosys} The
following properties hold:
\begin{enumerate}[$(a)$]
\item The functor $S$ preserves injective comodules.
\item If $N$ is a quasi-finite indecomposable right $eCe$-comodule then $S(N)$ is
indecomposable.
\item The functor $S$ preserves indecomposable injective comodules.
\item If $S_x$ is a simple $eCe$-comodule then $\Soc S(S_x)=S_x$.
\item If $S_x$ is a simple $eCe$-comodule then $S(S_x)$ is torsion-free.
\item We have that $S(\overline{E}_x)=E_x$ for all $x\in I_e$.
\item The functor $S$ preserves quasi-finite comodules.
\item The functor $S:\mathcal{M}^{eCe}\rightarrow \mathcal{M}^C$ restricts to
a fully faithful functor $S:\mathcal{M}^{eCe}_{qf}\rightarrow
\mathcal{M}^C_{qf}$ between the categories of quasi-finite comodules
which preserves indecomposables comodules and respects isomorphism
classes.
\end{enumerate}
\end{proposition}
\begin{proof}
\begin{enumerate}[$(a)$]
\item The functor $T$ is exact and left adjoint of $S$ so, by \cite[Proposition
9.5]{stenstrom}, the result follows.
\item Since $N$ is quasi-finite and indecomposable then we have that $\End_{eCe}(N)\cong \End_{C}(S(N))$ is a
local ring. Thus $S(N)$ is indecomposable.
\item It follows from $(a)$ and $(b)$.
\item Suppose that $\Soc S(S_x)=(\bigoplus_{i\in I} S_i)\bigoplus (\bigoplus_{j\in
J} T_j)$, where $S_i$ and $T_j$ are simple right $C$-comodules such
that $T(S_i)=S_i$ and $T(T_j)=0$ for all $i\in I$ and $j\in J$. By
Lemma \ref{soclequotient}, $\bigoplus_{i\in I} S_i=T(\Soc
S(S_x))\subseteq \Soc TS(S_x)=\Soc S_x=S_x$. Since $S$ is left exact
and preserves indecomposable injective comodules, $S_x\subseteq \Soc
S(S_x)\subseteq \Soc S(\overline{E}_x)=S_y$ for some simple comodule
$S_y$. Then $S_y=S_x=\Soc S(S_x)$.
\item If $M\subseteq S(S_x)$ is a non-zero torsion subcomodule of $S(S_x)$ then there exists a
simple $C$-comodule $R$ contained in $M$ such that $T(R)=0$. But
$\Soc S(S_x)=S_x$, so $S_x=R$ and we get a contradiction.
\item It is easy to see from $(c)$ and $(d)$.
\item Let $M$ be a quasi-finite right $eCe$-comodule. The injective envelope of $M$ is a
quasi-finite injective comodule $\xymatrix@1{M
\ar@<-0.3ex>@{^{(}->}[r] &E=\bigoplus \overline{E}_x^{n_x}}$. Since
$S$ is left exact then $\xymatrix@1{S(M) \ar@<-0.3ex>@{^{(}->}[r]
&S(E)=\bigoplus E_x^{n_x}}$. Thus $S(M)$ is quasi-finite.
\item It is a consequence of the above assertions and the equality $TS=1_{\M^{eCe}}$.
\end{enumerate}
\end{proof}

\begin{corollary}
$S$ preserves injective envelopes.
\end{corollary}

After proving Proposition \ref{injectivosys}, one should ask if the
behavior of simple comodules is analogous to injective ones, that
is, if $S$ preserves simple comodules and, consequently, in view of
Proposition \ref{injectivosys}$(c)$, $S(S_x)=S_x$ for all $x\in
I_e$. Unfortunately, in general, this is not true and we can only
say that $S(S_x)$ is a subcomodule of $E_x$ which contains $S_x$.

\begin{example}
This example shows that $S(S_x)$ does not have to be $S_x$ for every
$x\in I_e$. Consider the quiver $Q$
$$\xy \xymatrix@C=30pt{ \circ \ar[r]^-{\alpha} &
\circ} \POS (0,-2)*+{\txt{\scriptsize $y$}} \POS
(14.5,-2)*+{\txt{\scriptsize $x$}} \POS (17,0)*+{\txt{\scriptsize
,}}
\endxy$$
$C=KQ$ and the idempotent $e\in C^*$ associated to the set $\{x\}$.
Then, the localized coalgebra $eCe$ is $S_x$ and
$$S(S_x)=S_x\square_{eCe} Ce= eCe\square_{eCe}Ce\cong Ce\cong
<x,\alpha>\neq S_x.$$
\end{example}

The reader should observe that $S(S_x)$ could be an infinite
dimensional right $C$-comodule. Therefore, in general, $S$ cannot be
restricted to a functor between the categories of finite dimensional
comodules.

\begin{example} Consider the quiver $Q$
$$\xy \xymatrix@C=25pt@R=20pt{ & \cdots \ar[r]^-{\alpha_{n+1}}
&\circ \ar[r]^-{\alpha_{n}}& \circ \ar[r]^-{\alpha_{n-1}} & \circ
-\cdots - \circ \ar[r]^-{\alpha_2}&\circ  \ar[r]^-{\alpha_1} & \circ
 }
\POS (94.5,-3)*+{\txt{\scriptsize $1$}}
\POS (81.5,-3)*+{\txt{\scriptsize $2$}}
\POS (68.7,-3)*+{\txt{\scriptsize $3$}}
\POS (52.5,-3)*+{\txt{\scriptsize $n-1$}}
\POS (40,-3)*+{\txt{\scriptsize $n$}}
\POS (27,-3)*+{\txt{\scriptsize $n+1$}}
\POS (97,0)*+{\txt{\scriptsize ,}}
\endxy $$
$C=KQ$ and the idempotent $e\in C^*$ associated to the set $\{1\}$.
Then the localized coalgebra $eCe$ is $S_1$ and
$$S(S_1)=S_1\square_{eCe} Ce= eCe\square_{eCe}Ce\cong Ce\cong <1,\{
\alpha_1\cdots\alpha_{n-1}\alpha_n\}_{n\geq 1}>.$$

\end{example}

\begin{remark} The reader may find in \cite{jmn2} a proof of the
following fact: $S$ preserves finite dimensional comodules if and
only $S(S_x)$ is finite dimensional for each $x\in I_e$.
\end{remark}

In order to characterize the simple comodules invariant under the
functor $S$ we need the following result. It asserts that the
torsion immediate predecessors of a torsion-free vertex $S_x$ in
$\Gamma_C$ are the simple $C$-comodules contained in the socle of
$S(S_x)/S_x$. In the following picture the torsion-free vertices are
represented by white points.

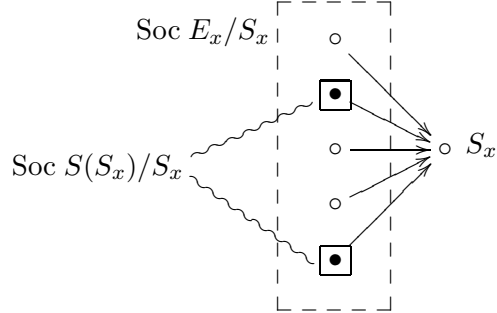
\[
  \xy \xymatrix@R=10pt@C=30pt{  & \circ  \ar[rdd]  &   \\
   & *+[F]{\bullet} \ar[rd] &  \\
     &   \circ \ar[r] &   \circ\\
      & \circ   \ar[ru]    &  \\
     &      *+[F]{\bullet}  \ar[ruu] &
   }

\POS (21,5) \ar@{--} (21,-36)

\POS (6,5) \ar@{--} (6,-36)

\POS (6,-36) \ar@{--} (21,-36)

\POS (6,5) \ar@{--} (21,5)

\POS (-6,-16) \ar@{~} (11,-8)

\POS (-6,-18) \ar@{~} (11,-30)

 \POS (33,-14)*+{\text{\small $S_x$}}

\POS (-4,1)*+{\text{\small $\Soc E_x/S_x$}}

\POS (-18,-17)*+{\text{\small $\Soc S(S_x)/S_x$}}

\endxy  \]

\begin{theorem}\label{torsionpredec}
Let $S_y$ and $S_x$ be two simple $C$-comodules. Then we have that
$S_y \subseteq S(S_x)/S_x$ if and only if $S_y \subseteq E_x/S_x$
and $T(S_y)=0$.
\end{theorem}
\begin{proof}
Consider the short exact sequence
\begin{equation}\label{equation}\xymatrix{S_x \ar[r] & S(S_x)
\ar[r]& S(S_x)/S_x }
\end{equation}
Since $S_x=T(S_x)=TS(S_x)$, $S(S_x)/S_x$ is a torsion subcomodule of
$E_x/S_x$. Therefore if $S_y \subseteq S(S_x)/S_x$ then $S_y
\subseteq E_x/S_x$ and $T(S_y)=0$.

Conversely, applying the functor $S$ to the exact sequence
$$\xymatrix{S_x \ar[r]^-i &
\overline{E}_x \ar[r]^-p&  \overline{E}_x/S_x}$$  we obtain the
following commutative diagram:
$$\xymatrix{S(S_x) \ar[r]^-{S(i)} &
E_x \ar[rr]^-{S(p)}\ar[rd]_-{\Coker S(i)}&  & S(\overline{E}_x/S_x)
\ar[r]&
\Coker\\
 &           &     E_x/S(S_x) \ar@<-0.3ex>@{^{(}->}[ru] \ar[r] & 0 &
}$$ Therefore we have that $\Hom_C(S_y,E_x/S(S_x))$ is contained in
the set of morphisms
$\Hom_C(S_y,S(\overline{E}_x/S_x))\cong\Hom_{eCe}(T(S_y),\overline{E}_x/S_x)=0$.
Now, applying $\Hom_C(S_y,-)$ to the exact sequence
$$\xymatrix{S(S_x) \ar[r]& E_x
\ar[r]& E_x/S(S_x),}$$ we obtain the exactness of the sequence
$$0=\Hom_C(S_y,E_x)\rightarrow \Hom_C(S_y,E_x/S(S_x)) \rightarrow
\Ext_C^1(S_y, S(S_x))\rightarrow 0$$
 and then
$0=\Hom_C(S_y,E_x/S(S_x))\cong \Ext_C^1(S_y,S(S_x))$.

Let us now apply the functor $\Hom_C(S_y,-)$ to (\ref{equation}) and
therefore $\Hom_C(S_y, S(S_x)/S_x)\cong \Ext_C^1(S_y, S_x)\neq 0$.
Then the result follows.
\end{proof}

\begin{corollary}\label{simples}
Let $S_x$ be a simple $eCe$-comodule. The following conditions are
equivalent:
\begin{enumerate}[$(a)$]
\item $E_x/S_x$ is torsion-free.
\item There is no arrow in $\Gamma_C$ from a torsion vertex $S_y$ to $S_x$.
\item $S(S_x)=S_x$.
\end{enumerate}
\end{corollary}
\begin{proof} It is straightforward from Lemma \ref{predecessors}
and Theorem \ref{torsionpredec}.
\end{proof}

We may generalize the former results using the Loewy series of the
$C$-comodule $S(S_x)$.

\begin{lemma}\label{techlemma} Let $S_x$ and $S_y$ be two simple $C$-comodules such
that $S_x$ is torsion-free. Then $S_y\subseteq
S(S_x)/\Socle^nS(S_x)$ if and only if $S_y$ is torsion and
$S_y\subseteq E_x/\Socle^nS(S_x)$.
\end{lemma}
\begin{proof} The necessity can be proved as the former
theorem. Conversely, since $S_y$ is torsion, by the proof of Theorem
\ref{torsionpredec}, we have that $\Ext^1_C(S_y,S(S_x))=0$. Then,
applying the functor $\Hom_C(S_y,-)$ to the short exact sequences
$$\xymatrix{\Socle^nS(S_x) \ar[r] & E_x \ar[r] &
E_x/\Socle^nS(S_x)}$$
$$\xymatrix{\Socle^nS(S_x) \ar[r] & S(S_x) \ar[r] &
S(S_x)/\Socle^nS(S_x)}$$ we obtain that
$\Hom_C(S_y,\frac{S(S_x)}{\Socle^nS(S_x)})\cong
\Hom_C(S_y,\frac{E_x}{\Socle^nS(S_x)})$.
\end{proof}

\begin{theorem}\label{propotorsionpred} Let $E_x$ be an indecomposable injective
$C$-comodule such that $S_x$ is torsion-free. If $S_y \subseteq
S(S_x)/\Socle^nS(S_x)$ for some $n\geq 1$, then the following
assertions hold:
\begin{enumerate}[$(a)$]
\item $S_y$ is torsion.
\item $S_y$ is a $n$-predecessor of $S_x$.
\item There exists a path in
$\Gamma_C$
%$$S_y\rightarrow S_{n-1}\rightarrow \cdots S_2\rightarrow
%S_1\rightarrow S_x$$
$$\xymatrix{S_y \ar[r] & S_{n-1} \ar[r] & \cdots \ar[r]&
S_2 \ar[r] & S_1 \ar[r] & S_x}$$
 such that $S_i$ is torsion for all $i=1,\ldots ,
n-1$.
\end{enumerate}
The converse also holds if $C$ is hereditary.
\end{theorem}
\begin{proof} $(a)$ can be proved as in Theorem \ref{torsionpredec}. $(b)$ is obtained from the inclusion
$\frac{S(S_x)}{\Socle^nS(S_x)}\subseteq \frac{E_x}{\Socle^n E_x}$.
For $(c)$, we proceed by induction on the number $n$. The case $n=1$
corresponds to Theorem \ref{torsionpredec}. Assume now that the
statement holds for $n-1$ and that $S_y \subseteq
S(S_x)/\Socle^nS(S_x)$. Analogously to the proof of Theorem
\ref{predecpath}, we may prove that there exists an arrow from $S_y$
to a simple $C$-comodule $S_j\subseteq S(S_x)/\Socle^{n-1}S(S_x)$.

Let us now suppose that $C$ is hereditary. We will prove that, for
each positive integer $n$, the torsion simple comodules contained in
$E_x/\Socle^n S(S_x)$ are those for which there is a path as
described in $(c)$. By Lemma \ref{techlemma}, this will imply the
statement. The case $n=1$ is just Theorem \ref{torsionpredec}. Let
us assume that it is verified for $n-1$. Then
$$\frac{E_x}{\Socle^{n} S(S_x)}\cong
\frac{\frac{E_x}{\Socle^{n-1}S(S_x)}}{\frac{\Socle^{n}S(S_x)}{\Socle^{n-1}S(S_x)}}\cong
\frac{\left (\oplus_{j\in J} E_j\right )\oplus \oplus
E}{\oplus_{j\in J} S_j}\cong \left (\oplus_{j\in J}E_j/S_j\right
)\oplus E,$$ where $S_j$ is torsion for all $j\in J$ and $E$ is an
injective comodule whose socle is torsion-free. If there is a path
from $S_y$ to $S_x$ as described in $(c)$, $1\in J$ by hypothesis
and hence $S_y\subseteq E_1/S_1\subseteq E_x/\Socle^n S(S_x)$.
Conversely, there is some $j_o\in J$ such that $S_y\subseteq
E_{j_0}/S_{j_0}$ and then we have an arrow $S_y\rightarrow S_{j_0}$.
By hypothesis there is a path of length $n-1$ as described in $(c)$
from $S_{j_0}$ to $S_x$. This completes the proof.
\end{proof}

\begin{corollary} Let $Q$ be a quiver, $KQ$ the path coalgebra
of $Q$ and $e\in (KQ)^*$ an idempotent element associated to the
subset $X\subseteq Q_0$. For each vertex $x\in X$, the $KQ$-comodule
$S(S_x)$ is generated by the set of paths
$$\xy
 \xymatrix@C=15pt{ \circ \ar[r] &  \circ \ar[r] &  \circ -\cdots
- \circ \ar[r] &  \circ \ar[r] & \circ }
\POS (0,-3)*+{\txt{\small $x_1$}}
\POS (10,-3)*+{\txt{\small $x_2$}}
\POS (43,-3)*+{\txt{\small $x_{n}$}}
\POS (54,-3)*+{\txt{\small $x$}}
\endxy $$
such that $x_i\notin X$ for any $i=1,\ldots, n$
\end{corollary}

\begin{corollary} Suppose that $\dim_K \Ext_C^1(S_1,S_2)$ is finite
for each simple $C$-comodules $S_1$ and $S_2$. If there is finitely
many paths in $\Gamma_C$, $$\xymatrix{S_t \ar[r] & S_{t-1} \ar[r] &
\cdots \ar[r]& S_2 \ar[r] & S_1 \ar[r] & S_x}$$ such that $S_i$ is a
torsion simple $C$-comodule for all $i=1,\ldots ,t$, then $S(S_x)$
is finite dimensional. If $C$ is hereditary, the converse holds.
\end{corollary}

\begin{remark} If we consider the non-valued Gabriel quiver of
$C$, it is not needed to assume that the groups of extensions
between simple comodules have finite dimension.
\end{remark}

\section{The quotient functor}\label{quotient}

Let us now analyze the properties of the quotient functor. We start
with an example which shows that, in general, $T$ does not preserve
injective comodules.

\begin{example}\label{example41}
Let $Q$ be the quiver
$$\xy\xymatrix@C=30pt{\circ \ar[r]^-{\alpha} & \circ
\ar[r]^-{\beta} & \circ }
\POS (0,-2)*+{\txt{\scriptsize $x$}}
\POS (14.5,-2)*+{\txt{\scriptsize $y$}}
\POS (29,-2)*+{\txt{\scriptsize $z$}}
\POS (32,0)*+{\txt{\scriptsize ,}}
\endxy$$
$C$ be the subcoalgebra of $KQ$ generated by $\{x,y,z, \alpha ,
\beta\}$ and $I_e=\{x,y\}$. The injective right $C$-comodule $E_z$
is generated by $<z,\beta>$ and $T(E_z)=<\beta>\cong S_y\neq E_y$.
\end{example}

\begin{proposition} The following statements hold:
\begin{enumerate}[$(a)$]
\item $T(E_x)=\overline{E}_x$ for any $x\in I_e$.
\item If $E$ is an injective torsion-free right $C$-comodule then
$T(E)$ is an injective right $eCe$-comodule.
\item If $M$ is a torsion-free right $C$-comodule then $\Soc M=\Soc T(M)=T(\Soc M)$.
\item The functor $T:\mathcal{M}^{C}\rightarrow \mathcal{M}^{eCe}$ restricts to
a functor $T:\mathcal{M}^{C}_{qf}\rightarrow \mathcal{M}^{eCe}_{qf}$
and a functor $T:\mathcal{M}^{C}_{f}\rightarrow
\mathcal{M}^{eCe}_{f}$ between the categories of quasi-finite and
finite dimensional comodules, respectively.
\end{enumerate}
\end{proposition}
\begin{proof}
 \begin{enumerate}[$(a)$]
\item By Proposition \ref{injectivosys}, $E_x=S(\overline{E}_x)$ for any $x\in I_e$. Then
$T(E_x)=TS(\overline{E}_x)=\overline{E}_x$.
\item It follows from $(a)$.
\item Consider the chain $\bigoplus_{x\in I} S_x=\Soc M\subseteq M\subseteq
E(M)=\bigoplus_{x\in I} E_x$. Since $M$ is torsion-free then
$I\subseteq I_e$. Therefore $\Soc M=\bigoplus_{x\in I}
S_x=\bigoplus_{x\in I} T(S_x)=T(\Soc M)\subseteq T(M)\subseteq
T(E(M))=\bigoplus_{x\in I}T(E_x)=\bigoplus_{x\in I} \overline{E}_x$
and the result follows.
\item It is easy to see.
\end{enumerate}
\end{proof}

\begin{corollary}
Let $E_x$ be a indecomposable injective $C$-comodule such that $S_x$
is torsion-free. $S(S_x)=E_x$ if and only if all predecessors of
$S_x$ in $\Gamma_C$ are torsion.
\end{corollary}
\begin{proof}
Assume that all predecessors of $S_x$ are torsion. Then,
$$T(\Socle^n E_x)=T(\Socle^{n+1} E_x)=T(\Soc E_x)=S_x.$$ Now, by the
previous proposition, $$\overline{E}_x=T(E_x)=T(\limdir \Socle^n
E_x)=\limdir T(\Socle^n E_x)=S_x$$ and thus, by Proposition
\ref{injectivosys}$(f)$, $E_x=S(\overline{E_x})=S(S_x)$.
 The converse follows from Theorem \ref{propotorsionpred}.
\end{proof}

\begin{example}\label{example410} In general, the functor $T$ is not full. Let $Q$ be the quiver
$$\xy \xymatrix@C=30pt{ \circ \ar[r]^-{\alpha} &
\circ} \POS (0,-2)*+{\txt{\scriptsize $y$}} \POS
(14.5,-2)*+{\txt{\scriptsize $x$}} \POS (17,0)*+{\txt{\scriptsize
,}}
\endxy$$
$C=KQ$ and $e\in C^*$ be the idempotent associated to the set
$\{x\}$. Then $\dim_K\Hom_C(S_x,C)=\dim_K\End(S_x)=1$ and
$\dim_K\Hom_{eCe}(S_x,eC)=2$. Therefore the map
$T_{S_x,C}:\Hom_C(S_x,C)\rightarrow \Hom_{eCe}(S_x, eC)$ cannot be
surjective.
\end{example}
\begin{example} In general, the functor $T$ does not preserve indecomposable
comodules. Let $KQ$ be the path coalgebra of the quiver
$$\xy \xymatrix@R=5pt@C=30pt{ \circ \ar[rd]^-{\alpha} & \\
& \circ\\
\circ \ar[ru]_-{\beta} & }
\POS (-2,0)*+{\txt{\scriptsize $x$}}
\POS (15,-8)*+{\txt{\scriptsize $z$}}
\POS (-2,-12)*+{\txt{\scriptsize $y$}}
\endxy$$
and $e\in C^*$ be the idempotent associated to the set $\{x,y\}$.
Then $T$ maps the indecomposable injective right $C$-comodule
$E_z=<z,\alpha, \beta>$ to the right $eCe$-comodule $S_x\oplus S_y$.
Nevertheless, it is easy to see that $T$ preserves indecomposable
torsion-free comodules.
\end{example}

Since $T(S_y)=0$ for each torsion simple $C$-comodule, one could
expect the analogous property for their injective envelopes.
Unfortunately, this is only true in some special conditions related
to the stability of the the torsion theory.

\begin{example} Let $KQ$ be the path coalgebra of the quiver
$$\xy \xymatrix@C=30pt{ \circ \ar[r]^-{\alpha} &
\circ,} \POS (0,-3)*+{\txt{\scriptsize $x$}} \POS
(15,-3)*+{\txt{\scriptsize $y$}}
\endxy$$
$C=KQ$ and $e\in C^*$ be the idempotent associated to the set
$\{x\}$. Then $T(E_y)=T(<y,\alpha>)\cong S_x\neq 0$.
\end{example}

\begin{theorem}\label{predec}
Let $E_y$ be an indecomposable injective right $C$-comodule with
$y\notin I_e$. The following statements are equivalent:
\begin{enumerate}[$(a)$]
\item $T(E_y)=0$,
\item $\Hom_C(E_y,E_x)= 0$ for all $x\in I_e$,
\item $S_y$ has no torsion-free predecessor in $\Gamma_C$.
\end{enumerate}
\end{theorem}
\begin{proof}
$(a) \Rightarrow (b).$ Since $S$ is left adjoint to $T$ then we have
that $\Hom_C(E_y,E_x)= \Hom_C(E_y,S(\overline{E}_x))\cong
\Hom_{eCe}(T(E_y), \overline{E}_x)=0$ for all $x\in I_e$.

$(b) \Rightarrow (c).$  It is proved in Proposition \ref{predechom}.

$(c) \Rightarrow (a).$ For each $n\geq 1$, we have the short exact
sequence
$$\xymatrix{\Socle^{n}E_y \ar[r]& \Socle^{n+1}E_y \ar[r] & \bigoplus_{i\in
I_n} S_i}.$$ Since $S_y$ has no torsion-free predecessor, $T(S_i)=0$
for all $i\in I_n$ and then $T(\Socle^{n}E_y)= T(\Socle^{n+1}E_y)$.
Now,
$$T(E_y)=T(\limdir
\Socle^nE_y)=\limdir T(\Socle^nE_y)=T(\Soc E_y)=0.$$
\end{proof}

Let us finish the section by giving an approach to the image of an
indecomposable injective comodule $E_y$ with torsion socle. Firstly,
from the Loewy series of $E_y$, we may obtain a chain
$$0\subseteq T(\Socle E_y)\subseteq T(\Socle^2 E_y)\subseteq\cdots
\subseteq T(\Socle^n E_y)\subseteq \cdots \subseteq T(E_y)$$ such
that each quotient $\frac{T(\Socle^{n+1} E_y)}{T(\Socle^{n}
E_y)}\cong T\left (\frac{\Socle^{n+1} E_y}{\Socle^{n} E_y}\right )$
is the the direct sum of the torsion-free $n$-predecessors of $S_y$.
As a consequent, by Lemma \ref{lemmasocle}, we have that
$T(\Socle^{n+1} E_y)\subseteq \Socle^n T(E_y)$.
 In particular, $T(\Socle^2 S_y)$ is
the direct sum of all torsion-free immediate predecessors of $S_y$
and $T(\Socle^2 E_y)\subseteq \Soc T(E_y)$.

\begin{lemma}\label{tdetorsion} Let $S_y$ be a torsion simple right $C$-comodule. Suppose
that $\{S_x, T_z\}_{x\in I, z\in J}$ is the set of all immediate
predecessors of $S_y$ in $\Gamma_C$, where $S_x$ is torsion-free for
all $x\in I$ and $T_z$ is torsion for all $z\in J$. Then
$$\Soc T(E_y)\subseteq \left (\bigoplus_{x\in I}
S_x\right ) \bigoplus \left (\bigoplus_{z\in J} \Soc T(E_z)\right
).$$ If $C$ is hereditary, the opposite inclusion also holds.
\end{lemma}
\begin{proof}
By Lemma \ref{predecessors}, $\Soc (E_y/S_y)=(\oplus_{x\in I}
S_x)\oplus (\oplus_{z\in J} T_z)$ and, consequently,
$E_y/S_y\subseteq (\oplus_{x\in I} E_x) \oplus (\oplus_{z\in J}
E_z)$. Then $$T(E_y)\cong T(E_y/S_y)\subseteq \left (\bigoplus_{x\in
I} \overline{E}_x\right ) \bigoplus \left (\bigoplus_{z\in J}
T(E_z)\right )$$ and then
$$\Soc T(E_y)\subseteq \left (\bigoplus_{x\in I} S_x\right ) \bigoplus \left (\bigoplus_{z\in J}
\Soc T(E_z)\right ).$$ Clearly, the inclusions are equalities if $C$
is hereditary.
\end{proof}

In a general context it is not possible to prove the equality in
Lemma \ref{tdetorsion}. For example, consider the quiver of Example
\ref{example41}, the coalgebra generated by the set
$<x,y,z,\alpha,\beta>$ and $I_e=\{x\}$. Then $\Soc T(E_z)=0$ and
$\Soc T(E_y)=S_x$.

\begin{corollary}\label{cort} Let $E_y$ be a indecomposable injective $C$
comodule such that $S_y$ is torsion. If $S_x\subseteq \Soc T(E_y)$
then \begin{enumerate}[$(a)$] \item $S_x$ is torsion-free. \item
$S_x$ is a predecessor of $S_y$ in $\Gamma_C$. \item There exists a
path in $\Gamma_C$
$$\xymatrix{S_x \ar[r] & S_{n} \ar[r] & \cdots \ar[r]&
S_2 \ar[r] & S_1 \ar[r] & S_y}$$ such that $S_i$ is torsion for all
$i=1,\ldots ,n$. \end{enumerate}

If $C$ is hereditary, the converse also holds.
\end{corollary}
\begin{proof} By Lemma \ref{tdetorsion}, it is enough to prove
$(b)$. Now, by hypothesis, $0\neq
\Hom_{eCe}(T(E_y),\overline{E}_x)\cong
\Hom_C(E_y,E_x)=\Rad_C(E_y,E_x)$. Corollary \ref{predechom}
completes the proof.
\end{proof}

\section{The colocalizing functor}\label{colocalizing}

Throughout this section we shall assume that $\T_e$ is a
colocalizing subcategory of $\M^C$. Then the quotient functor $T$
has a left adjoint functor $H$ which can be described as follows:

\[ \xymatrix@C=50pt{ \mathcal{M}^{C}
\ar@<0.75ex>[rr]^-{T=e(-)=-\square_{C}eC} & &
\ar@<0.75ex>[ll]^-{H=\Cohom_{eCe}(eC,-)} \mathcal{M}^{eCe}}.\]

We recall from \cite{takeuchi} that there exists such functor $H$ if
and only if $eC=T(C)$ is quasi-finite as right $eCe$-comodule, i.e.,
since $T(C)=eCe\bigoplus (\bigoplus_{y\in I_C\backslash I_e}
T(E_y))$, if and only if $\bigoplus_{y\in I_C\backslash I_e} T(E_y)$
is quasi-finite as $eCe$-comodule. According to Corollary
\ref{cort}, this is obtained if $\dim_K \Ext^1_C(S,S')$ is finite
for each pair of simple comodules $S$ and $S'$, and there is
finitely many paths in $\Gamma_C$
$$\xymatrix{S_x \ar[r] & S_{n} \ar[r] & \cdots \ar[r]&
S_3 \ar[r] & S_2 \ar[r] & S_1}$$ where $S_i$ is torsion for all
$i=1,\ldots ,n$, for each torsion-free simple comodule $S_x$.

\begin{proposition}\label{colocalizationfunctor} The following assertions hold:
\begin{enumerate}[$(a)$]
\item $H$ preserves projective comodules.
\item $H$ preserves finite dimensional comodules.
\item $H$ preserves finite dimensional indecomposable comodules.
\item The functor $H:\mathcal{M}^{eCe}\rightarrow \mathcal{M}^C$ restricts to
a fully faithful functor $H:\mathcal{M}^{eCe}_{f}\rightarrow
\mathcal{M}^C_{f}$ between the categories of finite-dimensional
comodules which preserves indecomposable comodules and respects
isomorphism classes.
\end{enumerate}
\end{proposition}
\begin{proof}
\begin{enumerate}[$(a)$]
\item  It is symmetric to the proof of Proposition
\ref{injectivosys}$(a)$.
\item Let $N$ be a finite dimensional right $eCe$-comodule. Then
$H(N)=\Cohom_{eCe}(eC,N)=\limdir \Hom_{eCe}(N_\lambda,
eC)^*=\Hom_{eCe}(N,eC)^*$. Now, since $eC$ is a quasi-finite right
$eCe$-comodule, $\Hom_{eCe}(N, eC)$ has finite dimension.
\item Let $N$ be a finite dimensional indecomposable right
$eCe$-comodule. Since $H$ is fully faithful then $\End_{eCe}(N)\cong
\End_{C}(H(N))$ is a local ring. Now, by $(b)$, $H(N)$ is finite
dimensional and then $H(N)$ is indecomposable.
\item It is straightforward from $(b)$, $(c)$ and the equality $TH=1_{\M^{eCe}}$.
\end{enumerate}
\end{proof}

Analogously to the study of the section functor, let us characterize
the simple comodules which are invariant under the functor $H$. For
that purpose we need the following proposition:

\begin{proposition}\label{lemmaprecoloc} Let $S_x$ be a simple $eCe$-comodule. Then $H(S_x)=S_x$ if
and only if $\Hom_{eCe}(S_x, T(E_y))=0$ for all $y\notin I_e$.
\end{proposition}
\begin{proof}
From the decomposition $eC=eCe\oplus eC(1-e)$ as right
$eCe$-comodules, we have the following equalities:
 $$\begin{array}{rl} \dim_K H(S_x)= &
\dim_K \Cohom_{eCe}(eC,S_x)\\ = &\dim_K S_x+\dim_K
\Cohom_{eCe}(eC(1-e),S_x)\\=& \dim_K S_x+\sum_{y\notin I_e} \dim_K
\Cohom_{eCe}(T(E_y),S_x)\\=& \dim_K S_x+\sum_{y\notin I_e}\dim_K
\Hom_{eCe}(S_x,T(E_y)). \end{array}$$ Therefore,
$\Hom_{eCe}(S_x,T(E_y))=0$ for all $y\notin I_e$ if $H(S_x)\cong
S_x$.

Conversely, there is a natural isomorphism $\Hom_C(H(S_x),S_x)\cong
\Hom_{eCe}(S_x,S_x)$ so there exists a non-zero (and then
surjective) morphism $f:H(S_x)\rightarrow S_x$. By hypothesis,
$dim_K H(S_x)=\dim_K S_x$ and then $f$ is an isomorphism.
%$$\begin{array}{rl} H(S_x)=S_x & \Leftrightarrow \Cohom_{eCe}(eC(1-e),S_x)=0 \\
% & \Leftrightarrow \Cohom_{eCe}(\oplus_{y\notin I_e} T(E_y),S_x)=0 \\
% & \Leftrightarrow \bigoplus_{y\notin I_e} \Cohom_{eCe}(T(E_y),S_x)=0  \\
% &\Leftrightarrow \bigoplus_{y\notin I_e}\Hom_{eCe}(S_x,T(E_y))=0 \\
% &\Leftrightarrow \Hom_{eCe}(S_x,T(E_y))=0 \hspace{0.2cm} \text{for all}\hspace{0.2cm}
%  y\notin I_e.\end{array}$$
\end{proof}

\begin{theorem}\label{propfunctorcoloc}
Let $S_x$ be a simple $eCe$-comodule. $H(S_x)=S_x$ if and only if
$\Ext^1_C(S_x,S_y)=0$ for all $y\notin I_e$, i.e., there is no arrow
$S_x\rightarrow S_y$ in $\Gamma_C$, where $S_y$ is a torsion simple
$C$-comodule.
\end{theorem}
\begin{proof}
By Proposition \ref{lemmaprecoloc}, it is enough to prove that
$\Ext^1_C(S_x,S_y)=0$ for all $y\notin I_e$ if and only if
$\Hom_{eCe}(S_x,T(E_y))=0$ for all $y\notin I_e$.

$\Leftarrow )$ Suppose that $\Ext^1_C(S_x,S_y)\neq 0$ for some
$y\notin I_e$. By Lemma \ref{predecessors}, $S_x\subseteq \Soc
(E_y/S_y)$. Then $$S_x=T(S_x)\subseteq T(\Soc E_y/S_y)\subseteq \Soc
T(E_y/S_y)=\Soc T(E_y)$$ and therefore $\Hom_{eCe}(S_x,T(E_y))\neq
0$.

$\Rightarrow )$ If $\Ext^1_C(S_x,S_y)= 0$ for each $y\notin I_e$ then there is no
 path as described in Corollary \ref{cort}$(c)$. Thus $S_x$ is not
contained in $T(E_y)$ for any $y\notin I_e$.

\end{proof}

\section{Semicentral idempotents}\label{secsemicentral}

In this section we extend the result obtained in \cite{jmnr} which
characterize stable subcategories as those whose associated
idempotent is left semicentral. We recall from \cite{jmnr} that an
idempotent element $e\in C^*$ is said to be \emph{left (right)
semicentral} if $eC=eCe$ ($Ce=eCe$), or equivalently, if $eC$ ($Ce$)
is a subcoalgebra of $C$. A localizing subcategory is said to be
\emph{stable} if it is closed for injective envelopes.

In the following theorem we describe stable subcategories from
different points of view. A proof of some equivalences is given in
\cite[Theorem 4.3]{jmnr}. We recall that, for a subset $\Lambda$ of
the vertex set $(\Gamma_C)_0$, we say that $\Lambda$ is \emph{right
link-closed} if it satisfies that, for each arrow $S\rightarrow T$
in $\Gamma_C$, if $S\in\Lambda$ then $T\in\Lambda$.

\begin{theorem}\label{stables} Let $C$ be a coalgebra and $\T_e\subseteq \M^C$ be a localizing
subcategory associated to an idempotent element $e\in C^*$. The
following conditions are equivalent:
\begin{enumerate}[$(a)$]
\item $\T_e$ is a stable subcategory.
\item $T(E_x)=0$ for any $x\notin I_e$.
\item $T(E_x)=\left\{%
\begin{array}{ll}
    \overline{E}_x & \text{if $x\in I_e$}, \\
    0 & \text{if $x\notin I_e$}. \\
\end{array}%
\right.$
\item $\Hom_C(E_y,E_x)=0$ for all $x\in I_e$ and $y\notin I_e$.
\item Any torsion vertex in $\Gamma_C$ has no torsion-free
predecessor.
\item $\mathcal{K}=\{S \in (\Gamma_C)_0 \mid T(S)=S\}$ is a right link-closed subset of $(\Gamma_C)_0$, i.e.,
 there is no arrow $S_x\rightarrow S_y$ in $\Gamma_C$, where $T(S_x)=S_x$ and $T(S_y)=0$.
\item There is no path in $\Gamma_C$ from a vertex $S_x$ to a vertex $S_y$ such that
$T(S_x)=S_x$ and $T(S_y)=0$.
\item $e$ is a left semicentral idempotent in $C^*$.
\end{enumerate}
\noindent If $T_e$ is a colocalizing subcategory these are also
equivalent to
\begin{enumerate}[$(i)$]
\item $H(S_x)=S_x$ for any $x\in I_e$.
\end{enumerate}
\end{theorem}
\begin{proof}
$(a)\Rightarrow (b)\Rightarrow (c)\Rightarrow (d)\Rightarrow
(e)\Rightarrow (f) \Rightarrow (g)$. It follows from the definition,
Theorem \ref{predec} and Proposition \ref{predechom}.

\noindent $(g) \Rightarrow (f)$. Trivial.

\noindent $(f) \Rightarrow (c)$. By Corollary \ref{cort}, $\Soc
T(E_y)=0$ for all $y\notin I_e$. Then $T(E_y)=0$ for all $y\notin
I_e$.

\noindent $(c)\Rightarrow (a)$. Let $M$ be a torsion right
$C$-comodule such that its injective envelope is $\bigoplus_{i\in
J}E_i$. Then $S_i\subseteq M$ is torsion for all $i\in J$ and, by
hypothesis, $T(E_i)=0$ for all $i\in J$. Thus $T(\bigoplus_{i\in
J}E_i)=0$.

\noindent $(c) \Leftrightarrow (h)$. We have that $C=\bigoplus_{x\in
I_C} E_x$ and therefore $T(C)=( \bigoplus_{x\in I_e} \overline{E}_x)
\bigoplus (\bigoplus_{y\notin I_e} T(E_y))$. On the other hand,
$eCe=\bigoplus_{x\in I_e} \overline{E}_x$. Therefore if $(c)$ holds
then $eCe=T(C)$. Conversely, if $eCe=T(C)$ then $\bigoplus_{x\in
I_e} \overline{E}_x=(\bigoplus_{x\in I_e} \overline{E}_x )\bigoplus
(\bigoplus_{y\notin I_e} T(E_y))$. Since $eCe$ is quasi-finite, by
Krull-Remak-Schmidt-Azumaya theorem, $T(E_y)=0$ for all $y\notin
I_e$.

\noindent $(e) \Leftrightarrow (i)$. It is Theorem
\ref{propfunctorcoloc}.
\end{proof}

We may find an analogous result to Theorem \ref{stables} for right
semicentral idempotents:

\begin{theorem} \label{rightsemicentral} Let $C$ be a coalgebra
 and $\T_e\subseteq \M^C$ be a localizing
subcategory associated to an idempotent element $e\in C^*$. The
following conditions are equivalent:
\begin{enumerate}[$(a)$]
\item $\T_{1-e}$ is a stable subcategory.
\item $T(E_x)=E_x$ for any $x\in I_e$.
\item $\Hom_C(E_x,E_y)=0$ for all $x\in I_e$ and $y\notin I_e$.
\item Any torsion-free vertex in $\Gamma_C$ has no torsion
predecessor.
\item There is no path in $\Gamma_C$ from a vertex $S_y$ to
a vertex $S_x$ such that $T(S_x)=S_x$ and $T(S_y)=0$.
\item $\mathcal{K}=\{S \in (\Gamma_C)_0 \mid T(S)=S\}$ is a left link-closed subset of
$(\Gamma_C)_0$, i.e., there is no arrow $S_y\rightarrow S_x$ in
$\Gamma_C$, where $T(S_x)=S_x$ and $T(S_y)=0$.
\item $e$ is a right semicentral idempotent in $C^*$.
\item The torsion subcomodule of a right $C$-comodule $M$ is $(1-e)M$
\item $S(S_x)=S_x$ for all $x\in I_e$.
\end{enumerate}
\end{theorem}
\begin{proof}
By \cite{jmnr} and Theorem \ref{stables}, it is easy to prove
$(a)\Leftrightarrow (b)\Leftrightarrow (c)\Leftrightarrow
(d)\Leftrightarrow (e) \Leftrightarrow (f)$.

$(d)\Leftrightarrow (g)$. It is Corollary \ref{simples}.
\end{proof}

\vspace{0.3cm} As a consequence, for a non necessarily
indecomposable coalgebra $C$, we get the following immediate result:
\begin{corollary}
The following are equivalent:
\begin{enumerate}[$(a)$]
\item $e$ is a central idempotent.
\item For each arrow $S_x\rightarrow S_y$ in $\Gamma_C$, $T(S_x)=0$ if and only if
$T(S_y)=0$.
\item For each connected component of $\Gamma_C$, either all vertices are torsion or all
vertices are torsion-free.
\item $T(E_x)=0$ for any $x\notin I_e$ and $S(S_x)=S_x$ for all $x\in I_e$.
\end{enumerate}
\end{corollary}

\begin{corollary}
Let $e$ be a central idempotent in $C^*$ and $\{\Sigma^t_C,
\Delta^s_C\}_{t,s}$ be the connected components of $\Gamma_C$, where
the vertices of each $\Sigma^t_C$ are torsion-free and the vertices
of each $\Delta^s_C$ are torsion. Then $\{\Sigma^t_C\}_t$ are the
connected components of $\Gamma_{eCe}$.
\end{corollary}

\end{document}